\documentclass[11pt]{article}

    \usepackage{charter}
    \usepackage[margin=1.2in]{geometry}
    \usepackage{amsmath, amsthm}
    \newtheorem{theorem}{Theorem}[section]  %
    \newtheorem{corollary}{Corollary}[theorem]  %

    \newtheorem{example}{Example}[section]

    \author{\small Gabriele Dragotto$^{1,3}$ \orcidA{} Stefan Clarke$^{1}$ \orcidB{} Jaime Fernández Fisac$^{2,3}$ \orcidC{} \\ \small Bartolomeo Stellato$^{1,2,3}$ \orcidD{} \\
    }
    \date{
    \emph{\small Department of Operations Research and Financial Engineering$^1$} \\
    \emph{\small Department of Electrical and Computer Engineering$^2$} \\
    \emph{\small Center for Statistics and Machine Learning$^3$} \\
    \emph{\small Princeton University}}

\usepackage[T1]{fontenc}
\usepackage[utf8]{inputenc}
\usepackage[usenames,dvipsnames]{xcolor}
\usepackage[colorlinks,
            linkcolor=BlueViolet,
            citecolor=blue,
            urlcolor=blue,
            bookmarks]{hyperref}
\usepackage{graphicx, pgfplots, subcaption, standalone, pgf, tikz}
\usepackage[justification=centering]{caption}
\usepackage[sort&compress,numbers,sectionbib]{natbib}
\usepackage{comment, url,float, setspace, xspace}
\pgfplotsset{compat=1.17}
\usepackage{array, bm, amssymb, amsfonts,  mathtools, nccmath}
\usepackage[inline]{enumitem}

\setlist[enumerate]{label=(\roman*.)}
\usepackage[ruled,vlined,linesnumbered]{algorithm2e}
\usepackage[nameinlink]{cleveref}
\usepackage{booktabs, longtable, multirow, multicol, adjustbox}

\definecolor{lime}{HTML}{A6CE39}
\DeclareRobustCommand{\orcidicon}{
    \hspace{-3mm}
	\begin{tikzpicture} 
    \draw[lime, fill=lime] (0,0) circle [radius=0.15] node[white] { 
        {\fontfamily{qag}\selectfont \tiny ID} 
    };
	\end{tikzpicture} 
    \hspace{-2mm}
}

\foreach \x in {A, ..., Z}{\expandafter\xdef\csname orcid\x\endcsname{\noexpand\href{https://orcid.org/\csname orcidauthor\x\endcsname}{\noexpand\orcidicon} } }

\newbool{boldSymbols}

\setbool{boldSymbols}{false}

\newcommand{\maxi}{\text{maximize}}
\newcommand{\mini}{\text{minimize}}
\newcommand{\st}{\text{subject to}}
\newcommand{\X}[1]{ \ifbool{boldSymbols}{\mathcal{X}}{X}^{ #1 }} 

\newcommand{\R}{\ifbool{boldSymbols}{\mathbb{R}}{\mathbf{R}}}
\newcommand{\Z}{\ifbool{boldSymbols}{\mathbb{Z}}{\mathbf{Z}}}
\newcommand{\Q}{\ifbool{boldSymbols}{\mathbb{Q}}{\mathbf{Q}}}
\newcommand{\NPhard}{\mbox{$\mathcal{NP}$-hard}\xspace}  
\newcommand{\identity}{\ifbool{boldSymbols}{\mathbf{I}}{I}}

\newcommand{\conv}{\operatorname{\bf conv}}

\newcommand{\ie}{\emph{i.e.}}
\newcommand{\eg}{\emph{e.g.}}

\crefname{line}{Line}{Lines}
\crefname{lemma}{Lemma}{Lemmata}
\crefname{theorem}{Theorem}{Theorems}
\crefname{proposition}{Proposition}{Propositions}
\crefname{algorithm}{Algorithm}{Algorithms}
\crefname{equation}{}{}
\crefname{definition}{Definition}{Definition}
\crefname{claim}{Claim}{Claim}
\crefname{corollary}{Corollary}{Corollaries}
\crefname{remark}{Remark}{Remarks}
\crefname{example}{Example}{Examples}
\crefname{figure}{Figure}{Figures}
\crefname{section}{Section}{Sections}
\crefname{table}{Table}{Tables}

\usepackage[xindy, acronyms, nowarn]{glossaries}   %
\glsdisablehyper
\makeglossaries
\newglossaryentry{MIP}
{
  name={MIP},
  plural={MIPs},
  description={mixed-integer optimzation problem},
  first={mixed-integer optimization problem (\glsentrytext{MIP})},
  descriptionplural={mixed-integer optimzation problem},
  firstplural={mixed-integer optimization problems (\glsentryplural{MIP})}
}
\newglossaryentry{CPL}
{
  name={CPL},
  plural={CPLs},
  description={cutting-plane layer},
  first={cutting-plane layer (\glsentrytext{CPL})},
  descriptionplural={cutting-plane layers},
  firstplural={cutting-plane layers (\glsentryplural{CPL})}
}

\newglossaryentry{CGP}
{
  name={CGP},
  plural={CGPs},
  description={cut-generating program},
  first={cut-generating program (\glsentrytext{CGP})},
  descriptionplural={cut-generating program},
  firstplural={cut-generating programs (\glsentryplural{CGP})}
}
\newglossaryentry{CGLP}
{
  name={CGLP},
  plural={CGLPs},
  description={cut-generating linear program},
  first={cut-generating linear program (\glsentrytext{CGLP})},
  descriptionplural={cut-generating linear programs},
  firstplural={cut-generating linear programs (\glsentryplural{CGLP})}
}
\newglossaryentry{RNN}
{
    name={RNN},
    description={recurrent neural network}
}

\newcommand{\parameters}{\ifbool{boldSymbols}{\boldsymbol{\theta}}{\theta}}
\newcommand{\cutparameters}{\ifbool{boldSymbols}{\boldsymbol{\sigma}}{\sigma}}
\newcommand{\rnnparameters}{\ifbool{boldSymbols}{\boldsymbol{\rho}}{\rho}}
\newcommand{\variables}{\ifbool{boldSymbols}{\mathbf{x}}{x}}

\newcommand{\integeridx}{\mathcal{I}}

\newcommand{\numconstraints}{m}
\newcommand{\numvars}{n}
\newcommand{\numtraining}{N}
\newcommand{\parameterset}{\Theta}
\newcommand{\probA}{\ifbool{boldSymbols}{\mathbf{A}}{A}(\parameters)}
  \newcommand{\probc}{\ifbool{boldSymbols}{\mathbf{c}}{c}(\parameters)}
  \newcommand{\probb}{\ifbool{boldSymbols}{\mathbf{b}}{b}(\parameters)}
\newcommand{\probregion}{X(\parameters)}
\newcommand{\probrelaxregion}{\tilde{X}(\parameters)}
\newcommand{\probexactregion}{\conv(\probregion)}
\newcommand{\optimal}{\ifbool{boldSymbols}{\mathbf{x}}{x}^\star(\parameters)}
\newcommand{\relaxsolution}[1]{\ifbool{boldSymbols}{\mathbf{x}}{x}^{#1}(\parameters)}
\newcommand{\toseparate}{\bar{\ifbool{boldSymbols}{\mathbf{x}}{x}}}
\newcommand{\relaxobjective}[1][]{z^{ #1 }(\parameters)}
\newcommand{\optimalobjective}{z^\star(\parameters)}
\newcommand{\cutsRHS}[1][]{\ifbool{boldSymbols}{\mathbf{h}}{h}^{#1}(\parameters)}
\newcommand{\cutsLHS}[1][]{\ifbool{boldSymbols}{\mathbf{G}}{G}^{#1}(\parameters)}
\newcommand{\cutg}{\ifbool{boldSymbols}{\mathbf{g}}{g}}
\newcommand{\cuth}{h}
\newcommand{\cutu}{\ifbool{boldSymbols}{\mathbf{u}}{u}}
\newcommand{\cutv}{\ifbool{boldSymbols}{\mathbf{v}}{v}}
\newcommand{\disj}{\ifbool{boldSymbols}{\boldsymbol{\pi}}{\pi}}
\newcommand{\disjvalue}{\ifbool{boldSymbols}{\pi_0}{\eta}}
\newcommand{\normalization}{\ifbool{boldSymbols}{\mathbf{D}}{D}}
\newcommand{\lr}{\mu}
\newcommand{\normtype}{p}

\newcommand{\rnnstate}[1]{\ifbool{boldSymbols}{\mathbf{S_{\rm state}}}{s_{\rm state}}^{#1}}
\newcommand{\algostate}[1][]{\ifbool{boldSymbols}{\mathbf{S}}{s}^{#1}(\parameters)}
\newcommand{\lenalgohistory}{M}
\newcommand{\rnndisj}{\ifbool{boldSymbols}{\boldsymbol{\alpha}}{\alpha}}
\newcommand{\rnnnorm}{\ifbool{boldSymbols}{\boldsymbol{\beta}}{\beta}}
\newcommand{\numcuts}{K}
\newcommand{\numiterations}{R}
\newcommand{\iteration}{r}
\newcommand{\loss}{L(\rnnparameters)}
\newcommand{\trainingset}{\mathcal{T}}
\newcommand{\ve}{\ifbool{boldSymbols}{v}{t}}
\newcommand{\uzero}{\ifbool{boldSymbols}{{u}_0}{u_{\numconstraints+1}}}
\newcommand{\vzero}{\ifbool{boldSymbols}{{v}_0}{v_{\numconstraints+1}}}
\newcommand{\subaddfn}{\Phi}
\newcommand{\doubleu}{\ifbool{boldSymbols}{\mathbf{w}}{w}}
\newcommand{\doubleuplus}{\ifbool{boldSymbols}{\mathbf{w}}{w}_+}
\newcommand{\doubleuminus}{\ifbool{boldSymbols}{\mathbf{w}}{w}_-}
\newcommand{\transpose}{\ifbool{boldSymbols}{\top}{T}}
\newcommand{\gomorydisjpart}{\ifbool{boldSymbols}{\mathbf{m}}{q}}
\newcommand{\lemmay}{\ifbool{boldSymbols}{\mathbf{y}}{y}}
\newcommand{\onezerovec}{\ifbool{boldSymbols}{\mathbf{r}}{\xi}}
\newcommand{\utilde}{\tilde{\ifbool{boldSymbols}{\mathbf{u}}{u}}}
\newcommand{\vtilde}{\ifbool{boldSymbols}{\mathbf{v}}{v}}
\newcommand{\spareMone}{\ifbool{boldSymbols}{\mathbf{m}}{q}_1}
\newcommand{\spareMtwo}{\ifbool{boldSymbols}{\mathbf{m}}{q}_2}
\newcommand{\spareMthree}{\ifbool{boldSymbols}{\mathbf{m}}{q}_3}
\newcommand{\spareMfour}{\ifbool{boldSymbols}{\mathbf{m}}{q}_4}
\newcommand{\spareMfive}{\ifbool{boldSymbols}{\mathbf{m}}{q}_5}

\newcommand{\fractionalpart}[1]{\{ #1 \}}
\newcommand{\zerovec}{\ifbool{boldSymbols}{\mathbf{0}}0}

\newcommand{\facilityvars}{\ifbool{boldSymbols}{\mathbf{y}}{y}}
\newcommand{\facilitycosts}{\ifbool{boldSymbols}{\mathbf{c}}{c}(\parameters)}

\newcommand{\controltime}{\tau}
\newcommand{\controlmaxtime}{H}
\newcommand{\power}{P}
\newcommand{\maxpower}{\power^{\text{max}}}
\newcommand{\samplingtime}{\upsilon}
\newcommand{\powerload}{\power^{\text{load}}}
\newcommand{\energy}{E}
\newcommand{\minenergy}{\energy^{\text{min}}}
\newcommand{\maxenergy}{\energy^{\text{max}}}
\newcommand{\initenergy}{\energy^{\text{init}}}
\newcommand{\fuelstate}{\zeta}

\newcommand{\betacontrol}{\omega}
\newcommand{\gammacontrol}{\delta}
\newcommand{\dcontrol}{\xi}
\newcommand{\hcontrol}{h}
\newcommand{\wcontrol}{\psi}
\newcommand{\fuelstateinit}{\fuelstate^{\text{init}}}
\newcommand{\scontrol}{k}
\newcommand{\scontrolinit}{s^{\text{init}}}

\newcommand{\nsw}{N^{\text{sw}}}
\newcommand{\Gcontrol}{G}
\newcommand{\controlf}{f}
\newcommand{\dpast}{\dcontrol^{\text{past}}}

\title{ Differentiable Cutting-plane Layers for Mixed-integer Linear Optimization }

\begin{document}
\maketitle

\begin{abstract}
We consider the problem of solving a family of parametric mixed-integer linear optimization problems where some entries in the input data change. We introduce the concept of \emph{cutting-plane layer} (CPL), \ie, a differentiable cutting-plane generator mapping the problem data and previous iterates to cutting planes. We propose a CPL implementation to generate split cuts, and by combining several CPLs, we devise a differentiable cutting-plane algorithm that exploits the repeated nature of parametric instances. In an offline phase, we train our algorithm by updating the internal parameters controlling the CPLs, thus altering cut generation. Once trained, our algorithm computes, with predictable execution times and a fixed number of cuts, solutions with low integrality gaps. Preliminary computational tests show that our algorithm generalizes on unseen instances and captures underlying parametric structures.
\end{abstract}

\section{Introduction}
We consider the problem of optimizing a \emph{parametric family} of \glspl{MIP} of the form
\begin{equation}
    \label{eq:mip}
     \optimalobjective = 
     \begin{array}[t]{ll}
     \mini  & \probc^\transpose \variables\\
     \st  & \probA \variables \le  \probb, \quad \variables \ge 0,\\ & \variables_i \in \Z,\quad \forall i \in \integeridx,
    \end{array}
\end{equation}
with respect to the decision vector $\variables \in \R^\numvars$, where the $\integeridx$ is the set of indices of integer-constrained variables, and the rational problem data $\probA \in \R^{\numconstraints \times \numvars}$, $\probb \in \R^{\numconstraints}$, and $\probc \in \R^\numvars$ depend on the parameters $\parameters \in \parameterset$. We denote the feasible region of problem~\eqref{eq:mip} as $\probregion$, and we assume that $\probA$, $\probb$ contain the upper bounds on $\variables$.
For any given $\parameters \in \parameterset$, problem~\cref{eq:mip} represents a single instance, which we assume to be feasible and bounded.

In the past decades, we have witnessed dramatic theoretical and computational developments in solving single \glspl{MIP} instances, in particular, linear ones~\cite{bixby2010}.
Researchers and practitioners can nowadays choose among several commercial~\cite{gurobi,fico,cplex} and open-source~\cite{scip, highs} mixed-integer solvers based on branch-and-cut~\cite{padberg_branch-and-cut_1991}, an expert blend of the branch-and-bound~\cite{land_automatic_1960} algorithm, cutting planes \cite{gomory_outline_1958,kelley1960cutting}, and heuristics. 
While the above methodologies are tremendously effective for solving single \glspl{MIP} from scratch, little is known, in theory and practice, on how to optimize a family of similar  \glspl{MIP}~\cite{mipcc23,warmstart}. 
In many real-world settings, users repeatedly solve similar \emph{parametric} instances where only some critical entries of the input data change, \eg, prices, demands, or capacities. For instance, an energy company may generate the optimal power dispatch by repeatedly solving a unit commitment problem where only the predictions for energy demand vary across time. 
Similarly, a hybrid (\eg, mixed-integer) controller may solve an optimal control problem to update the system inputs (\eg, actuator torques) as the sensor signals (\eg, system state) and goals (\eg, desired trajectory) vary over time.
Whenever the variables are all continuous, linear-programming duality enables us to understand how the parametric variations in $\parameters$ affect solutions; in principle, duality extends to the discrete case as \emph{subadditive duality}, yet the approach ``does not integrate well with current computational practice''~\cite{guzelsoy2007duality}.
Therefore, the majority of the limited literature on parametric \glspl{MIP} and reoptimization focuses on reusing the information of the branch-and-cut tree (\eg, \cite{gamrath2015reoptimization,patel2023progressively}); it comes as no surprise that exploiting this information poses significant unsolved mathematical and computational challenges in mixed-integer optimization~\cite{mipcc23} and beyond~\cite{amos_amortized_2023}.

\paragraph{Our contributions.} 
This work opens a novel research direction by introducing a differentiable cutting-plane algorithm to optimize parametric families of \glspl{MIP}. 
We summarize our contributions as follows:

\begin{enumerate}
\item We introduce a differentiable cutting-plane algorithm designed to exploit the repeated structure of parametric instances (\cref{sec:algorithm}). 
During an offline phase, we train our algorithm to generate effective cutting planes, precisely split cuts, by \emph{dynamically adjusting} the \emph{cut-generating parameters}, \ie, internal parameters controlling how cuts are generated (\eg, disjunction and normalization). 
Our algorithm interprets cutting-plane generation as a {\it sequential decision problem} where the {\it control inputs} are the cut-generating parameters, and the {\it state} is the accumulated information across cutting-plane rounds, \eg, the generated cuts and the primal bounds. 
Once trained, our algorithm can tackle unseen instances by computing, in a fixed number of cutting plane rounds, solutions with low integrality gaps.
\item We introduce \glspl{CPL}, \ie, differentiable cutting-plane generators with tunable weights (\cref{sec:layers}). 
\glspl{CPL} are the basic building blocks of our algorithm, and support a forward pass, \ie, the generation of cutting planes performed by solving a convex optimization problem, and a backward pass, \ie, differentiation with respect to the layers' weights.
We propose an implementation of \glspl{CPL} to generate split cuts, and we establish a theoretical correspondence between the split cuts and the generalized subadditive Gomory cuts from Ch\'etelat and Lodi~\cite{chetelat_continuous_2023}.
\item Finally, we conduct preliminary testing of our algorithm on two small families of parametric 2-matching and hybrid control instances where either the objective or right-hand side changes (\cref{sec:experiments}). Our algorithm generalizes well on unseen instances and provides up to a $3$-fold improvement over a standard cutting-plane algorithm.
\end{enumerate}
\subsection{A concise literature review}
\paragraph{Learning for mixed-integer optimization.}

A nurtured stream of research leverages machine learning to enhance mixed-integer optimization~\cite{MachineLearninBengio2021,dalle2022learning}, focusing mainly on data-driven heuristics for branch-and-cut algorithms (\eg, \cite{learn2branch,alvarez2017,Lodi2017,exactcombopt,guzelsoy2007duality,bertsimas2021,stellato2022,cauligi2021coco}). Several works~\cite{huang2022learning,faenza_cutting_2020,paulus2022learning,berthold2022learning,turner2023cutting,turner2022adaptive,deza_cutting_2023} addressed the problem improving cutting plane algorithms, primarily focusing on \emph{cutting plane selection}, \ie, the problem of selecting when and which cuts to add; see~\cite{deza_cutting_2023} for a detailed survey. 
While these heuristics offer appealing computational advantages, they are limited to selecting from a finite set of commonly used cuts.
In addition, they rely on a complete branch-and-cut algorithm, whose overhead can be prohibitively large for some applications, \eg, real-time embedded optimization. 
In this work, we propose a fixed-length cutting-plane algorithm that dynamically tailors cut generation to a specific family of parametric problems. 
In sharp opposition to the previous works, we do not perform cutting plane selection but \emph{generation},  \ie, we control the generation of the cuts and their coefficients.
To the best of our knowledge, the only exception is Ch\'etelat and Lodi~\cite{chetelat_continuous_2023}, where the authors propose a continuous cutting-plane algorithm based on a specific class of subadditive functions. 
While this approach shares several characteristics with ours (see \cref{sec:layers:theory}), it focuses on solving a single problem and is not designed for parametric families.

\paragraph{Generating good cuts.}
The classes, number, and general properties of the cuts significantly affect the efficiency and convergence of cutting-plane algorithms~\cite{padberg_branch-and-cut_1991,dey2018theoretical,cptutorial}. 
Although popular metrics for cut effectiveness include cut violation, parallelism, or sparsity~\cite{bixby2002solving,dey2018theoretical,walter2013sparsity,amaldi2014coordinated}, they primarily focus on the effect of the cuts at the current round, which may not be a direct metric for the overall algorithm behavior.
Our algorithm will, instead, learn to construct cuts by maximizing the bound improvement they induce across the entire sequence of cutting plane rounds (\eg, similarly to Padberg's \emph{optimal separator}~\cite[Ch.9]{padberg2013linear}).
In addition, we choose to focus on disjunctive cuts, whose popularity soared in correspondence with their implementation inside branch-and-cut frameworks~\cite{balas1980strengthening,balas_1993_lift,balas1996mixed,ontheseparation} and their conic extensions (\eg, \cite{kilincc2016minimal,lodi2023disjunctive}). 
Specifically, we consider \emph{split cuts}~\cite{splitcuts} because, first, they include many other cuts as special cases (\eg, Gomory mixed-integer and mixed-integer rounding cuts); and second, because their closure can tightly approximate the convex-hull of $\probregion$~ \cite{fischetti2013approximating,balas2008optimizing,dash2010mir}. As reported by Cornu\'ejols~\cite{cptutorial}, ``the split closure is surprisingly tight, and integer programming solvers should probably do a better job of approximating it''. 

\paragraph{Learning sequential decision-making.}
We model our cutting plane algorithm as a sequential decision problem where the cut-generating parameters represent the \emph{input actions}, and the continuous relaxation represents a \emph{dynamical system} (\eg, an environment), both interacting in a closed-loop fashion. When we cannot model the system or its goal with a mathematical model, reinforcement learning techniques, \eg, policy gradient methods, can improve the closed-loop behavior~\cite{sutton2018reinforcement,reinforce,fisac2019bridging}.
However, in our setting, the environment and the goal (\eg, increasing relaxation's bound) are well-defined, and structure-exploiting approaches offer more efficient alternatives in terms of data and computation~\cite{moerland2023model}.
In a sense, our \glspl{CPL} act as convex-optimization control policies~\cite{cocp} that can be efficiently trained and evaluated by solving small optimization problems.

\section{Differentiable cutting plane algorithm}
\label{sec:algorithm}

Consider a single instance from the parametric family of \glspl{MIP} in \cref{eq:mip}  by fixing~$\parameters$. We denote its (an) optimal solution as $\optimal$ with objective $\optimalobjective$, and we refer to $\probexactregion$ as the convex hull of its feasible solutions. 
We can solve the instance via a cutting-plane algorithm, an iterative method that progressively tightens a relaxation of problem~\cref{eq:mip} by adding extra inequalities. 
Specifically, in each round (iteration) $\iteration$, a standard cutting plane algorithm solves the \emph{relaxed} problem
\begin{equation}
    \label{eq:relax}
     \relaxobjective = \begin{array}[t]{ll}
     \mini & \probc^\transpose \variables\\
    \st &  \probA \variables \le  \probb,
    \quad \cutsLHS \variables \le \cutsRHS,\quad \variables \ge 0,
    \end{array}
\end{equation}
where we denote the feasible region as~$\probrelaxregion$.
In this linear optimization problem, we drop the integrality requirements, and we add a set of inequalities $\cutsLHS \variables \le \cutsRHS$,  where $\cutsLHS \in \R^{\numcuts \times \numvars}$ and $\cutsRHS \in \R^{\numcuts}$.
Such inequalities, called \emph{cutting planes} or \emph{cuts}, guarantee that at each round, the previous fractional solution is cut off while retaining any integer-feasible point in $\probexactregion$.
More precisely, given a solution $\relaxsolution{\iteration}$ to problem~\cref{eq:relax}, we say that an inequality $\cutg^\transpose \variables \le \cuth$, with $\cutg \in \R^{\numvars}$ and $\cuth \in \R$, is
\begin{enumerate*}
\item a valid inequality if it is satisfied by any point of~\eqref{eq:mip}, \ie, $\cutg^\transpose \variables \le \cuth$ for any $x \in \probexactregion$, and 
\item a cut for $\relaxsolution{\iteration}$ if it is valid for $\probexactregion$ and violated by $\relaxsolution{\iteration}$, \ie, $\cutg^\transpose \relaxsolution{\iteration} > \cuth$.
\end{enumerate*}
Since $\relaxobjective$ is a lower bound for $\optimalobjective$, cutting plane algorithms guarantee monotonic bound improvements.

\subsection{Our algorithm} 
Our cutting-plane algorithm should satisfy two foremost properties: it should provide \emph{high-quality solutions} for the parametric family of \glspl{MIP} in \cref{eq:mip}, and it should be \emph{differentiable}. 
On the one hand, when we say provide high-quality solutions, we mean that the algorithm should dynamically adapt the cut generation (\ie, separation) procedure to exploit the structure of the parametric family; 
for instance, by identifying which constraints or variables should play a more active role in cut generation.
On the other hand, to be \emph{differentiable}, our algorithm should implement both a forward pass, \ie, unrolling a given number of iterations of a cutting plane method, and a backward pass, \ie, computing gradients with respect to its key parameters. 
Similarly to \cite{chetelat_continuous_2023}, differentiability plays a pivotal role in adjusting the parameters that control cutting-plane generation; we call these \emph{cut-generating parameters} $\cutparameters$. 
\begin{figure}[!ht]
\includegraphics[width=\textwidth]{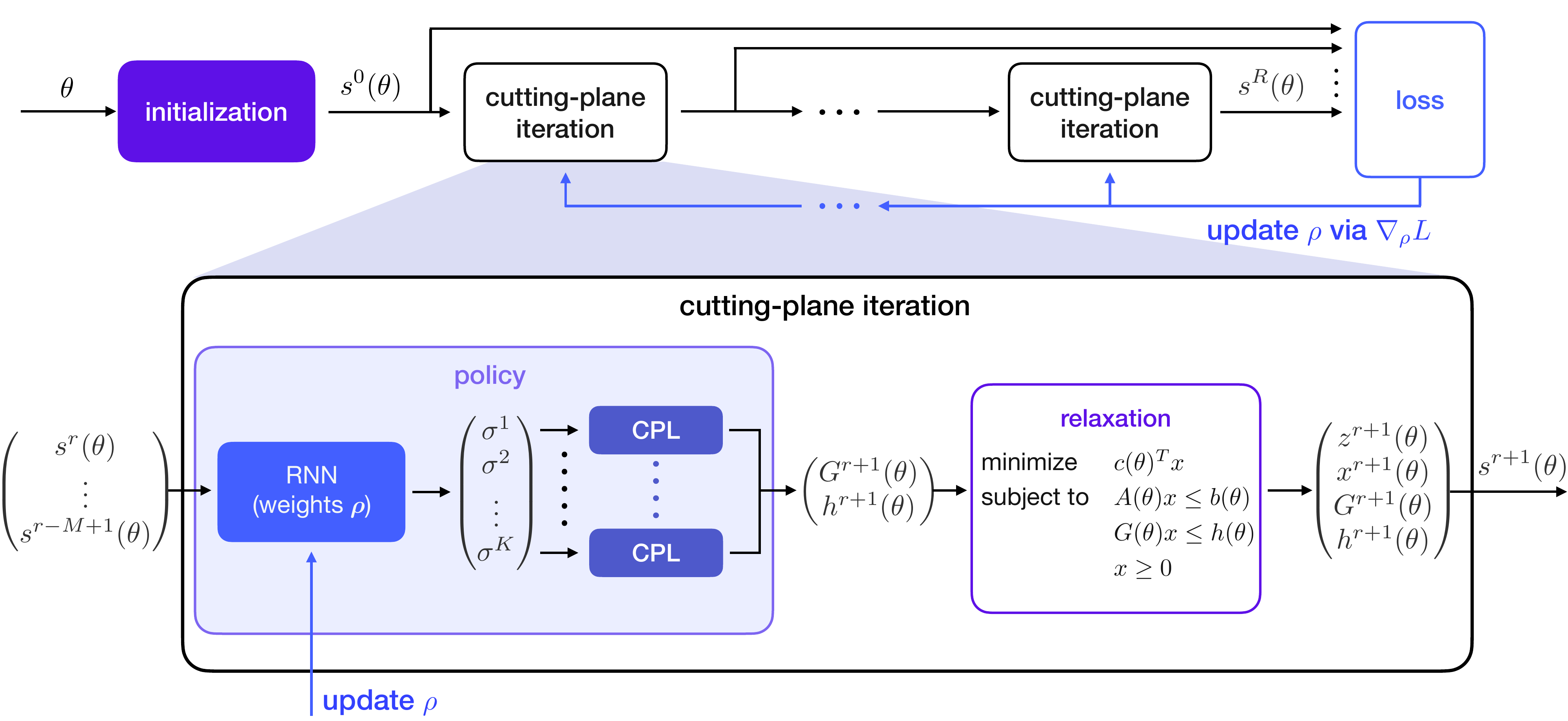}
\caption{Overview of our differentiable cutting-plane algorithm. }
\label{fig:architecture}
\end{figure}
\paragraph{Cutting plane algorithm.} 
Our algorithm, outlined  in~\cref{fig:architecture}, has a fixed size, meaning it performs a finite number of cutting plane rounds $\numiterations$ where, in each round, it adds $\numcuts$ cutting planes. 
As an input, the algorithm receives a vector~$\parameters$ defining a particular instance of problem~\eqref{eq:mip} and produces, as an output, a solution $\relaxsolution{\numiterations}$.
During each iteration $\iteration=0,\dots,\numiterations$,
we represent the state of our algorithm as $\algostate[\iteration] = (\relaxobjective[\iteration],\relaxsolution{\iteration}, \cutsLHS[\iteration], \cutsRHS[\iteration]),$ \ie, the current objective value, candidate solution, and cutting plane coefficients.
The first step, the {\it initialization block}, initializes $\algostate[0]$ by solving the initial continuous relaxation without cutting planes. Given this state, our algorithm generates $K$ cuts and solves the resulting continuous relaxation, obtaining a new state $\algostate[\iteration+1]$.

\paragraph{Policy.}
The idea of using a policy accounting for the \emph{state} of the cutting-plane algorithm stems from the observation that cut generation is, indeed, an optimal control problem~\cite{amaldi2014coordinated}; as previously generated cuts affect the generation in later iterations, the cut generator should exploit the state of the system to produce better cuts. 
Here, a \emph{policy} generates the cut coefficient (\ie, control inputs) using two components. First, an \gls{RNN} produces cut parameters $(\cutparameters^1,\dots,\cutparameters^K)$ from the history of past $\lenalgohistory$ states $(\algostate[\iteration], \dots,\algostate[\iteration - \lenalgohistory + 1])$.
Second, we feed the output of the \gls{RNN} into $K$ \glspl{CPL} to generate the cuts $(\cutsLHS[\iteration + 1], \cutsRHS[\iteration + 1]),$ as detailed in \cref{sec:layers}.
Our algorithm then solves the continuous relaxation to produce a new candidate solution $\relaxsolution{\iteration+1}$ with objective $\relaxobjective[\iteration+1]$, and an updated state $\algostate[\iteration+1]$.

\paragraph{Loss function.} 
To optimize the performance of our architecture for the considered family of parametric instances defined in problem~\eqref{eq:mip}, we need to train the policy such that it generates high-quality cutting plane coefficients, by tuning the \gls{RNN} weights, denoted as $\rnnparameters$.
Specifically, our goal is to minimize the \emph{loss function},
\begin{align}
\label{eq:loss}
\loss = {\bf E}_{\theta} 
\left(-\sum_{\iteration=1}^{\numiterations}  \gamma^r \probc^\transpose \left( \relaxsolution{\iteration} - \relaxsolution{\iteration-1}  \right) \right),
\end{align}%
which is the expectation over $\theta$ of the negative sum of primal bound improvements, discounted by factor~$\gamma \in (0, 1)$.
This loss function is {\it unsupervised}, as it does not include the true optimal solution but only the bound induced by the cuts. 
Intuitively, a higher primal bound in early iterations correlates with improved policy performance, leading to a reduced loss. 
Furthermore, the \gls{RNN} will likely produce cut parameters, and thus cuts, inducing greater improvements in the primal bound over $\numiterations$ rather than, for example, improvements in the violation with respect to the previous candidate solution.

\paragraph{Training.}
We assume to observe a set of the realization of the parameters $\parameters$, which we will call the \emph{training set} $\trainingset = \{\parameters_i\}_{i=1}^{\numtraining}$.
To train our policy, we use a stochastic gradient method~\cite{optmethodsml} of the form $\rnnparameters^{\ell + 1} = \rnnparameters^{\ell} - \lr \nabla_{\rnnparameters}\hat{L}(\rnnparameters^{\ell}),$ where $\nabla_{\rnnparameters}\hat{L}(\rnnparameters^{\ell})$ is an unbiased stochastic gradient of $\loss$ with respect to the \gls{RNN} weights $\rnnparameters$, and $\lr$ is the learning rate.
This process computes the gradients using the chain rule (\ie, backpropagation through time~\cite{werbos1990backpropagation}), which also requires differentiating through optimization problems (\eg, the continuous relaxation and the \glspl{CPL}); the latter differentiation is equivalent to the implicit differentiation of each optimization problem's optimality conditions~\cite{diffcp2019,cvxpylayers2019}.

\section{Cut generation and cutting plane layers}
\label{sec:layers}

We consider split cuts, \ie, disjunctive cuts from two-sided disjunctions, which are generated by exploiting a Farkas argument on the original simplex tableau or via the so-called \gls{CGLP} \cite{balas_disjunctive_2018,splitcuts}. We focus on the latter, as it enables us to explore a larger space of cutting planes, \eg, via expressive normalizations (see \cref{sec:layers:why}).
Since this problem is convex, we can differentiate its optimal solution with respect to its parameters \cite{diffcp2019,cvxpylayers2019}.

\paragraph{The CGP.}
We slightly generalize the \gls{CGLP} to a parametrized convex optimization problem in the cut parameters $\cutparameters$, which we denote as the \gls{CGP}. Compared to \gls{CGLP}, the \gls{CGP} features a convex normalization constraint, but the separation procedure is analogous to that of \gls{CGLP}. Specifically, we generate cuts by considering a \emph{valid} two-sided \emph{split} disjunction for $\conv(\probregion)$ in the form of $(\disj^\transpose \variables \le \disjvalue) \lor (\disj^\transpose \variables \ge \disjvalue +1)$ with $\disj \in \R^\numvars$ and $\disjvalue \in \R$. When we say the disjunction is valid, we mean that any feasible solution in $\probregion$ belongs to one side of the disjunction. At a given iteration $\iteration$, we can determine if there exists a cut $\cutg^\transpose \variables \le \cuth$ separating a fractional solution $\toseparate$ from 
$\conv (
    \{ \variables \in \probrelaxregion : \disj^\transpose \variables \le \disjvalue \}
    \cup 
    \{\variables \in \probrelaxregion : \disj^\transpose \variables \ge \disjvalue +1 \}
)$
by solving the \gls{CGP} formulated~as
\begin{subequations}
\label{eq:CGP}
	\begin{align}
        \underset{{\cutg,\cuth,\cutu,\cutv}}{\maxi\;} \quad &\cutg ^\transpose \toseparate - \cuth 
        \label{eq:CGP:objective} \\[-3pt]
		\st \quad & \cutg \le \begin{bmatrix}  \probA^\transpose & \disj  \end{bmatrix}\cutu, \\[-3pt]
        & \cutg \le \begin{bmatrix}   \probA^\transpose & -\disj \end{bmatrix}\cutv, \\[-3pt]
        & \cuth \ge \begin{bmatrix}  \probb^\transpose & \disjvalue \end{bmatrix}\cutu, \\[-3pt]
        & \cuth \ge \begin{bmatrix}  \probb^\transpose & -(\disjvalue+1)\end{bmatrix}\cutv, \\[-3pt]
        &\left\lVert \normalization\begin{bmatrix} \cutu \\ \cutv \end{bmatrix}\right\lVert_\normtype\le1, \label{eq:CGP:normalization}\\[-3pt]
        & \cutu \ge 0, \cutv \ge 0
	\end{align}
\end{subequations}
where $\cutu \in \R^{\numconstraints+1}$, $\cutv \in \R^{\numconstraints+1}$, $\cutg \in \R^{\numvars}$, $\cuth \in \R$, $\normalization \in \R^{(\numconstraints+1)\times(\numconstraints+1)}$, and $\normtype$ is the norm type. 
The cut parameters are $\cutparameters = (\disj, \disjvalue, \normalization)$.
The \gls{CGP} is a linear problem (the usual \gls{CGLP}) if $\normtype = 1$ or $\infty$, and a second-order cone problem if $\normtype=2$. 
By construction, any feasible solution $(\cutg, \cuth, \cutu, \cutv)$ to problem~\cref{eq:CGP} such that $\cutg ^\transpose \toseparate - \cuth>0$ corresponds to a violated split cut. 
Without \cref{eq:CGP:normalization}, the feasible set of \cref{eq:CGP} is a cone, and if there exists at least one solution such that $\cutg ^\transpose \toseparate - \cuth>0$, the problem is unbounded. Therefore, we truncate the cone via the normalization constraint~\cref{eq:CGP:normalization} to recover a cut. 
Traditionally, the normalization constraint involves a $0/1$ linear combination of $\cutu$ and $\cutv$ (or $\cutg, \cuth$), and has the role of selecting which cut to generate by pushing \gls{CGLP} to select one of the extreme rays. For instance, if $\normtype =1$ and $\normalization$ is a zero matrix except $\normalization_{\numconstraints + 1, \numconstraints + 1} = \normalization_{2(\numconstraints+1), 2(\numconstraints+1)} = 1$,
we recover the so-called \emph{trivial} normalization.
Instead, in \cref{eq:CGP:normalization}, we consider a generalized normalization constraint  involving the $\normtype$-norm of a linear transformation of $\cutu$ and $\cutv$ via $\normalization$. 
One natural question is whether generalizing the \gls{CGLP} via a 2-norm brings any benefits. 
When $\normtype=1$, cuts derived from the extreme rays of \gls{CGLP} dominate those in the interior \cite{ontheseparation}. 
Empirically, as we show in \cref{sec:experiments}, \emph{dominated cuts} may lead to better cuts in later iterations. 
Furthermore, as we motivate in the next section, $\normtype=2$ enables us to represent a large family of (subadditive) cuts.

\begin{comment}
\paragraph{Strengthening cuts.} As long as we can differentiate through elementary operations (\eg, the flooring function), we can strengthen $\cutg$ and $\cuth$. Let $\onezerovec =\{0\}^m \times \{1\}$, and  $\probA_j$ be the $j$-th column of $\probA$. We strengthen the coefficient $\cutg_j$ of each integer variable $j$ via a monoidal strengthening \cite{balas1980strengthening} as
%
\begin{align*}
    %
    \bar{\cutg}_j &= 
        \max \Bigg(\cutu^\transpose \begin{bmatrix}
        \probA_j \\ \lfloor \mathbf{m}_j \rfloor + \disj_j
    \end{bmatrix}, \cutv^\transpose \begin{bmatrix}
        \probA_j \\ -\lceil \mathbf{m}_j \rceil + \disj_j
    \end{bmatrix}\Bigg), & \mathbf{m} = \frac{(\cutu - \cutv)^\transpose}{(\mathbf{1} - \onezerovec)(\cutu + \cutv)} \begin{bmatrix}
        \probA_j \\ \disj_j
    \end{bmatrix}.
\end{align*}
\end{comment}

\subsection{Cut representation and range}
\label{sec:layers:theory}

\paragraph{Subadditive valid inequalities.}
For brevity, we limit our exposition to the pure-integer case, although our results extend to the mixed-integer one, and we focus on a single realization of~$\parameters$. Let the fractional part of a vector $\lemmay$ be  $\fractionalpart{\lemmay} = \lemmay - \lfloor \lemmay \rfloor$.
Ch\'etelat and Lodi \cite{chetelat_continuous_2023} devise a cut generation scheme via a \emph{subadditive} function, parametrized in $\doubleu \in \R^\numconstraints$ and $\ve \in \R$, defined as
\begin{align*}
    \subaddfn_{\doubleu, \ve}(\lemmay) &= \min \Big( \fractionalpart{\doubleu^\transpose \lemmay}, \frac{\ve}{1 - \ve}(1 - \fractionalpart{\doubleu^\transpose \lemmay})\Big) + \max \Big(-\doubleu, \frac{\ve}{1 - \ve} \doubleu \Big)^\transpose \lemmay.
\end{align*}
The resulting inequality $\subaddfn_{\doubleu, \ve}(-\probA) \variables \ge \subaddfn_{\doubleu, \ve}(-\probb)$ is valid for $\conv(\probregion)$ for any $\doubleu \in \R^\numconstraints$ and  $\ve \in [0, 1)$. 
In \cref{thm:continuouscp}, we show that the subadditive inequalities of \cite{chetelat_continuous_2023} are feasible solutions to problem~\cref{eq:CGP}, and thus split cuts with specific disjunctions.
In the following, we write $\probA_j$ to represent the $j$th column of $\probA$. 
We also write the sum of the function $\subaddfn_{\doubleu, \ve}$ over the \emph{columns} of $\probA$ as $\subaddfn_{\doubleu, \ve}(\probA) = \sum_{j=1}^\numconstraints \subaddfn_{\doubleu, \ve}(\probA_j)$.
\begin{theorem}
    \label{thm:continuouscp}
    Let $\doubleu \in \R^\numconstraints$ and $\ve \in [0, 1)$ be such that $\ve = \fractionalpart{\probb^\transpose \doubleu }$. Let $\doubleuplus \ge 0$ and $\doubleuminus \ge 0$ be the positive and negative parts of $\doubleu$, \ie, $\doubleu = \doubleuplus - \doubleuminus$ and ${\doubleuplus}^\transpose \doubleu_- = 0$. Let $\vzero = \ve$, $\uzero = 1 - \ve$, $\disjvalue = \lfloor\probb^\transpose \doubleu\rfloor$, and, for any variable~$j$,
    \begin{align*}
        \disj_j &= \begin{cases}
            \lfloor -\probA_j^\transpose \doubleu \rfloor &  \probA_j^\transpose \doubleu \le 0 \\
            \lceil - \probA_j^\transpose \doubleu \rceil & \probA_j^\transpose \doubleu > 0
        \end{cases}.
        \end{align*}
    Finally, let $\cutu = (\doubleuplus,  \uzero),$ $\cutv = (\doubleuminus, \vzero),$ $\cuth = [\probb^\transpose  -(\disjvalue+1)] \cutv$, and
    \begin{equation*}
        \cutg = \min\Big(\begin{bmatrix} \probA^\transpose &\disj \end{bmatrix} \cutu,\quad \begin{bmatrix} \probA^\transpose &-\disj \end{bmatrix} \cutv \Big).
    \end{equation*}
    If \cref{eq:CGP:normalization} represents the trivial normalization, then
    \begin{enumerate}
        \item \label{thm:continuouscp:1} $(\cutg, \cuth, \cutu, \cutv)$ is feasible in \cref{eq:CGP} and  defines the valid inequality $\cutg^\transpose \variables \le \cuth$,
        \item \label{thm:continuouscp:2} The inequality $\cutg^\transpose \variables \le \cuth$ is equivalent to $\subaddfn_{\doubleu, \ve}(-\probA) \variables \ge \subaddfn_{\doubleu, \ve}(-\probb)$,
        \item \label{thm:continuouscp:3} For any other $\hat{\ve} \in [0, 1) \neq \ve$, the inequality $\subaddfn_{\doubleu, \ve}(-\probA) \variables \ge \subaddfn_{\doubleu, \ve}(-\probb)$ dominates $\subaddfn_{\doubleu, \hat{\ve}}(-\probA) \variables \ge \subaddfn_{\doubleu, \hat{\ve}}(-\probb)$.
    \end{enumerate}
\end{theorem}
We refer to the appendix for a proof of \ref{thm:continuouscp:1} and \ref{thm:continuouscp:2}, whereas \ref{thm:continuouscp:3} follows from \cite{chetelat_continuous_2023}.
In view of these results, the subadditive inequalities from~\cite{chetelat_continuous_2023} can be interpreted as an \emph{explicit} \gls{CPL} in the cut-generating parameters~$\doubleu$ and~$\ve$, while \cref{eq:CGP} defines an \emph{implicit} \gls{CPL}. 
For implicit, we mean that cuts are the solutions to the optimization problem rather than the result of a closed-form expression. Naturally, there are advantages and disadvantages in either approach. The subadditive inequalities are \emph{unoptimized} cuts corresponding to feasible, yet not necessarily optimal, solutions of \cref{eq:CGP}. This implies that they are often dominated, weak, and not properly cuts but rather valid inequalities.
That is why the approach in~\cite{chetelat_continuous_2023} requires an over parametrization of their generating parameters, \ie, over parametrization of $\doubleu,\ve$ into a matrix and a vector, respectively, such that many, in the hundreds, valid inequalities are generated. 
Conversely, our explicit \gls{CPL} optimizes over the space of cuts and induces a monotonic bound in the differentiable cutting-plane algorithm; however, it requires us to solve a convex optimization problem to compute each cut.

\paragraph{Cut range.}  Since our approach generates arbitrary split cuts, it can also generate Gomory mixed-integer, lift-and-project, and mixed-integer rounding cuts \cite{cptutorial,nemhauser1990recursive,cornuejols2002rank}. 
Hence, \gls{CGP} with $\normtype=1$, even in the restricted space of diagonal non-negative matrices $\normalization$, can already represent substantial families of cuts. 
We provide an additional corollary related to the types of cuts we can represent with $\normtype=2$. 
We refer to the appendix for the complete proof.

\begin{corollary}
    Let $\toseparate$ be a fractional point. If $\normtype=2$ in \cref{eq:CGP:normalization}, then the set of cuts generated by varying $\cutparameters=(\normalization, \disj, \disjvalue)$ and solving \cref{eq:CGP} contains the set of undominated subadditive cuts cutting off $\toseparate$. 
    \label{coro:norm2}
\end{corollary}
The previous corollary implies that the closure of the cuts generated by \cref{eq:CGP} with $\normtype=2$ contains the closure of the subadditive Gomory cuts
$\subaddfn_{\doubleu, \ve}(-\probA) \variables \ge \subaddfn_{\doubleu, \ve}(-\probb)$ over  $\doubleu \in \R^m, \ve = \fractionalpart{\doubleu^\transpose \probb}$.

\subsection{Why cut-generating parameters matter}
\label{sec:layers:why}
We briefly explain the motivations, intuitions, and implications of dynamically adjusting~$\cutparameters$.

\paragraph{Disjunctions.} First, the policy can select the variables participating in the split disjunction via $\disj$ and $\disjvalue$, and thus in the generation of the cut.
The problem of selecting how to form the disjunctions in $\disj$ is generally \NPhard \cite{ralphs_selecting_2010,caprara2003separation}. Balas and Sexena \cite{balas2008optimizing} optimize over the space of possible split disjunctions by solving a parametric mixed-integer cut-generating program. Our split \glspl{CPL} can be seen as a learning version of \cite{balas2008optimizing} that produces a distilled number of cuts. However, in contrast with \cite{balas2008optimizing}, we do not solve mixed-integer problems as subroutines, and, importantly, we let the normalization $\normalization$ also vary as part of $\cutparameters$.

\paragraph{Normalization.} Different normalizations $\normalization$ push the \gls{CGP} to select different extreme rays and the corresponding basis in the relaxation \cref{eq:relax}. Indeed, \gls{CGP} with $\normtype=1$ can be interpreted as a method for finding a good basis from which to generate a Gomory cut \cite{bonami2012optimizing}. 
However, our expressive normalization induces a broader range of cuts. Our intuition is that, on the one hand, the cuts generated in the parametric family most likely derive from the same subset of constraints, corresponding to the active multipliers $\cutu, \cutv$ in \cref{eq:CGP}. Therefore, by adjusting the normalization coefficients, we practically change the \emph{cost} of using certain constraints in the cut-generation process. On the other hand, we generate cuts that, although seemingly weak, may result in stronger cuts in later rounds. As we show in \cref{ex:normalization}, this seems to be the case even in small examples.

\begin{example}
\label{ex:normalization}
Consider Example 1 from Fischetti et al.~\cite{lodi2023disjunctive}, a binary problem with $\numvars=2$. One of the possible split (lift-and-project) inequalities generated from the initial continuous relaxation is $x_1+4x_2 \le 1$. However, this cut is never optimal for \gls{CGLP} under traditional normalization conditions (\ie, standard or trivial normalization). However, we found it is optimal for several normalizations \cref{eq:CGP:normalization} where $\normtype=1$ and $\normalization$ is a diagonal matrix. 
In \cref{fig:normalization}, we train our differentiable cutting-plane algorithm to specifically solve this problem ($\numiterations=2$ with $\numcuts=1$ cut each). We represent the cutting planes as the light-blue regions. Although $x_1+4x_2 \le 1$ is not sufficient to solve the problem to optimality, a tailored, slowed-down version of our algorithm tilts the second cut in only three backward passes (1 backward pass is enough).
\begin{figure}[!ht]
\centering
\begin{subfigure}{0.4\textwidth}
  \centering
  \includegraphics[width=.8\linewidth]{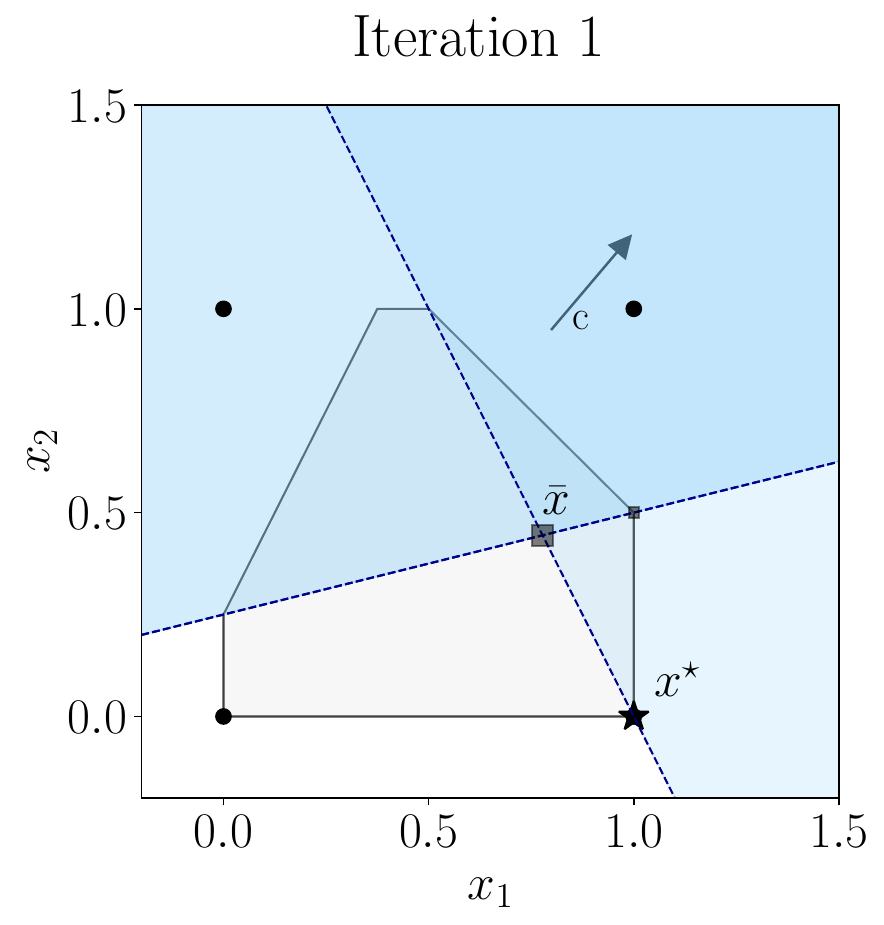}
  \label{fig:sfig1}
\end{subfigure}%
\begin{subfigure}{0.4\textwidth}
  \centering
  \includegraphics[width=.8\linewidth]{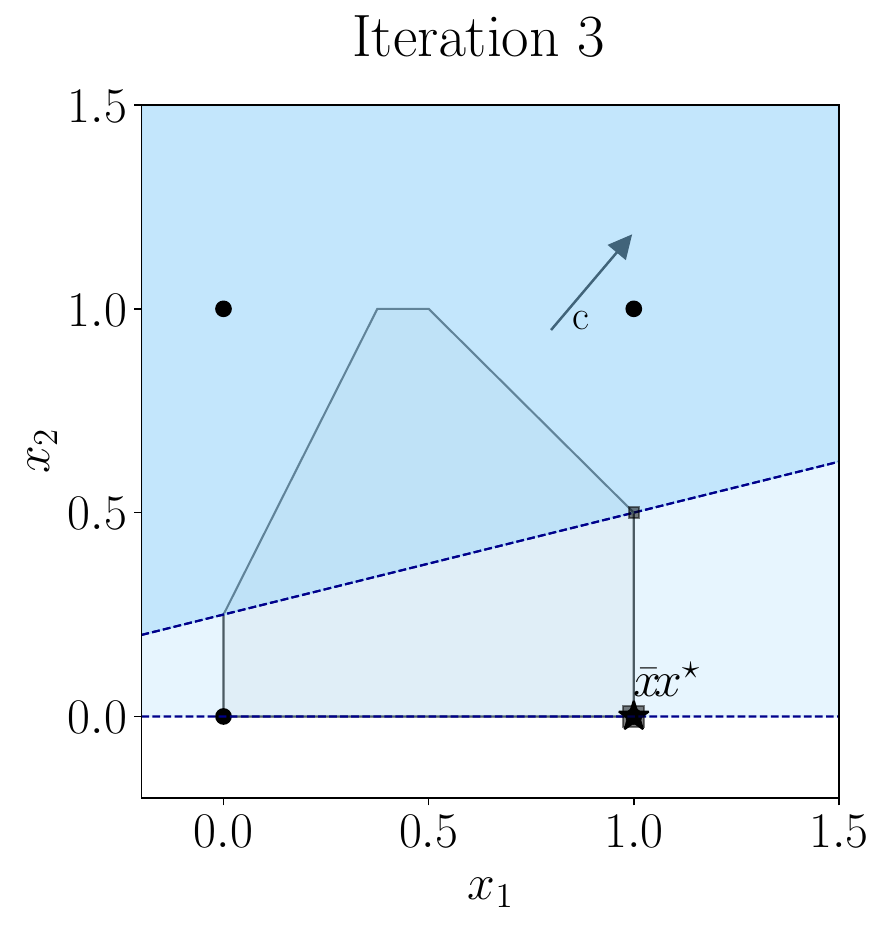}
  \label{fig:sfig2}
\end{subfigure}
\caption{Example 1 from Fischetti et al. \cite{lodi2023disjunctive} solved with 3 backward passes. }
\label{fig:normalization}
\end{figure}
\end{example}
With a $2$-norm, we gain some computational advantages in terms of differentiation as the norm regularizes the \gls{CGP} into a better-behaved second-order cone problem from a differentiation perspective. 
Furthermore, it generates cut coefficients that may increase cuts' diversity and thus propel greater improvements in the primal bound. 
We empirically validate this hypothesis in the next section.

\section{Computational experiments}
\label{sec:experiments}
We evaluate our differentiable cutting-plane algorithm on two parametric families of 2-matching and hybrid control problems. Although our tests rely on unoptimized and experimental libraries not designed for \gls{MIP} (\eg, with limited sparse linear algebra support), they provide promising numerical insights.

\paragraph{Experimental setup.} We employ stochastic gradient descent with momentum and a decaying learning rate to perform the backward pass, and we set the discount factor~$\gamma$ in \cref{eq:loss} to $0.9$. 
We implement our algorithm in \texttt{JAX}~\cite{jax2018github}, we use \texttt{Clarabel}~\cite{clarabel} (tolerance of $10^{-8}$) as an interior-point solver, and \texttt{cvxpylayers}~\cite{cvxpylayers2019} to differentiate \cref{eq:relax} and \cref{eq:CGP}. 
We run our tests on a single-core Intel Xeon Gold $6142$ with $32$GB of RAM. 
We define our \gls{RNN} as an LSTM cell layer of size $512$, followed by $2$ fully-connected layers of sizes $4(\numconstraints+1)+2\numvars$ and $2(\numconstraints+1)+\numvars$, and a softmax output activation function. 
We restrict our normalization coefficients $\normalization$ to diagonal matrices.
In each training \emph{epoch}, we iterate over the entire training set, which we split into batches of $10$ instances.
For any $\numcuts$, $\numiterations$, we compare our algorithm against the same number and rounds of split cuts generated with a \gls{CGLP} normalized via the standard normalization (\emph{baseline}). We strengthen the cuts, in both cases, via a monoidal strengthening \cite{balas1980strengthening}. We select our models based on the performance of our algorithm on a \emph{validation} set of instances, and we perform validation-based early-stop. Finally, we compare the baseline and our algorithm on a hidden \emph{test} set. We summarize, in \cref{tab:results}, our results on the training, validation, and test datasets. Specifically, we report the mean integrality gap (\emph{Gap}), mean integer infeasibilities (\emph{Infeas}), and maximum constraint violations (\emph{MaxViol}) with respect to the candidate solution in each round.

\subsection{2-Matching instances}
\begin{figure}[!ht]
\includegraphics[width=\textwidth]{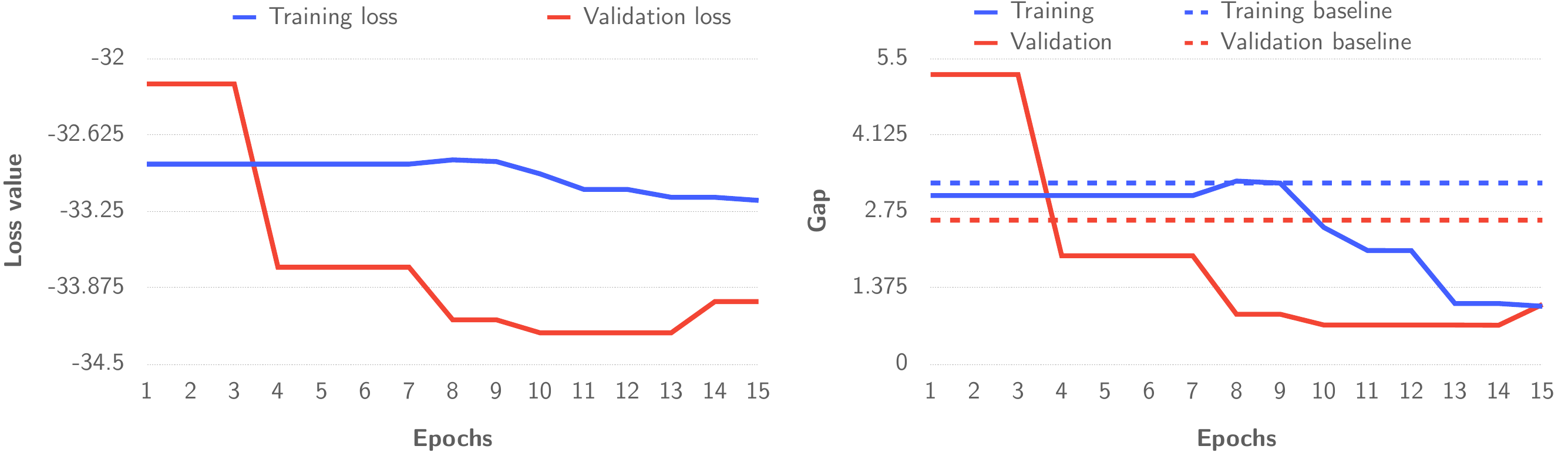}
\caption{Training and validation results for parametric 2-matchings.}
\label{fig:matching}
\end{figure}

We consider $70$ small 2-matching instances generated from a random Erdős–Rényi graph where $\numvars=35$ and $\numconstraints=16$ (excluding bounds). We sample $\parameters$ by varying the objective coefficient of each variable, corresponding to an edge, by sampling it over a normal distribution $\mathcal{N}(30,50)$ truncated at the third decimal. In these matching instances, Chvátal-Gomory cuts, and thus Gomory cuts, are sufficient to describe the convex hull of feasible solutions \cite{chvatal1973edmonds}.
The average integrality gap for the training, validation, and test datasets are $43.45\%$, $37.88\%$, and $38.73\%$, respectively, while their cardinality is $\numtraining=40$, $15$, and $15$ instances each, respectively.  We train our algorithm with a learning rate of $10^{-5}$, a norm $\normtype=2$, a disjunction of only one variable (\ie, lift-and-project), $\numiterations=3$, and $\numcuts=2$.
\vspace{-1em}
\begin{table}[ht]
\resizebox{\columnwidth}{!}{%
\begin{tabular}{l@{\hspace{1em}}l@{\hspace{2em}}r@{\hspace{1em}}r@{\hspace{1em}}r@{\hspace{2em}}r@{\hspace{1em}}r@{\hspace{1em}}r@{\hspace{2em}}r@{\hspace{1em}}r@{\hspace{1em}}r}
                           & & \multicolumn{3}{c}{\textbf{Training}}                      & \multicolumn{3}{c}{\textbf{Validation}}                    & \multicolumn{3}{c}{\textbf{Test}}                          \\
                           & & \textbf{\small Gap} & \textbf{\small Infeas} & \textbf{\small MaxViol} & \textbf{\small Gap} & \textbf{\small Infeas} & \textbf{\small MaxViol} & \textbf{\small Gap} & \textbf{\small Infeas} & \textbf{\small MaxViol} \\  \toprule
\multirow{2}{*}{Matchings} & \textbf{CPLs} & 1.09         & 0.01                   & 0.50               & 0.70         & 0.01                   & 0.49               & 0.00         & 0.00                   & 0.50               \\
                          & \textbf{Baseline}          & 3.26         & 0.01                   & 0.05               & 2.59         & 0.03                   & 0.49               & 1.17         & 0.02                   & 0.50   \\ \hline
\multirow{2}{*}{Control} & \textbf{CPLs} & 0.38     & 0.17          & 0.30      & 2.28       & 0.17          & 0.25      & 5.47  & 0.14          & 0.13      \\
                         & \textbf{Baseline}               & 18.10    & 0.05          & 0.16      & 21.60      & 0.05          & 0.19      & 21.65 & 0.05          & 0.08     \\
                         \bottomrule\\
\end{tabular}
}
\caption{Overview of computational results. \label{tab:results}}
\end{table}
\vspace{-3em}
\paragraph{Results.}  We summarize our results in \cref{fig:matching,tab:results}. The baseline gaps for training, validation, and tests are $3.26\%$, $2.59\%$, and $1.17\%$, respectively.  Our algorithm quickly improves the cuts and converges in about $15$ epochs despite starting from the same gap as the baseline.  We observe a training, validation, and test gap of $1.09\%$, $0.70\%$, and $0.00\%$, a $3$-fold improvement over the baseline. This set of tests also highlights the role of the norm-$2$ as a regularizer for \gls{CGP}. As speculated, we believe weaker cuts in earlier rounds lead the algorithm to find better cuts in later rounds. We empirically observed that the maximum violation of the candidate solution monotonically decreases as the training progresses. 
Therefore, cut violation, even in a well-known 2-matching setting, might not be the most appropriate proxy for evaluating cut effectiveness.

\subsection{Optimal control}
We consider some parametric hybrid control problems modeling a fuel cell and a capacitor system according to the formulations in \cite{stellato2022,frick2015embedded}, without the small quadratic objective term ($\le 10^{-3}$).
We refer to the appendix for the details on this model. 
Here, a set of binary variables control the on/off status of a fuel cell that balances the energy flow with a capacitor, and linear equations determine the system's dynamics. 
We model a control horizon of $10$, resulting in instances where $\numconstraints=90$ (excluding bounds) and $\numvars=40$. We vary the constraints' right-hand sides ($\parameters$) corresponding to the initial conditions of the system dynamics. 
We randomly sample $0/1$ values for the integer-constrained parameters, whereas, for the continuous parameters, we sample from a normal distribution centered in the average between their lower and upper bounds.
The average integrality gaps for the training, validation, and test datasets are $21.83\%$, $23.04\%$, and $24.17\%$, respectively, while their cardinality is $\numtraining=50$, $25$, and $25$, respectively. 
We train our algorithm with a learning rate of $10^{-5}$, $\normtype=1$, a disjunction of only one variable (\ie, lift-and-project), $\numiterations=2$ and $\numcuts=5$.

\begin{figure}[!ht]
\includegraphics[width=\textwidth]{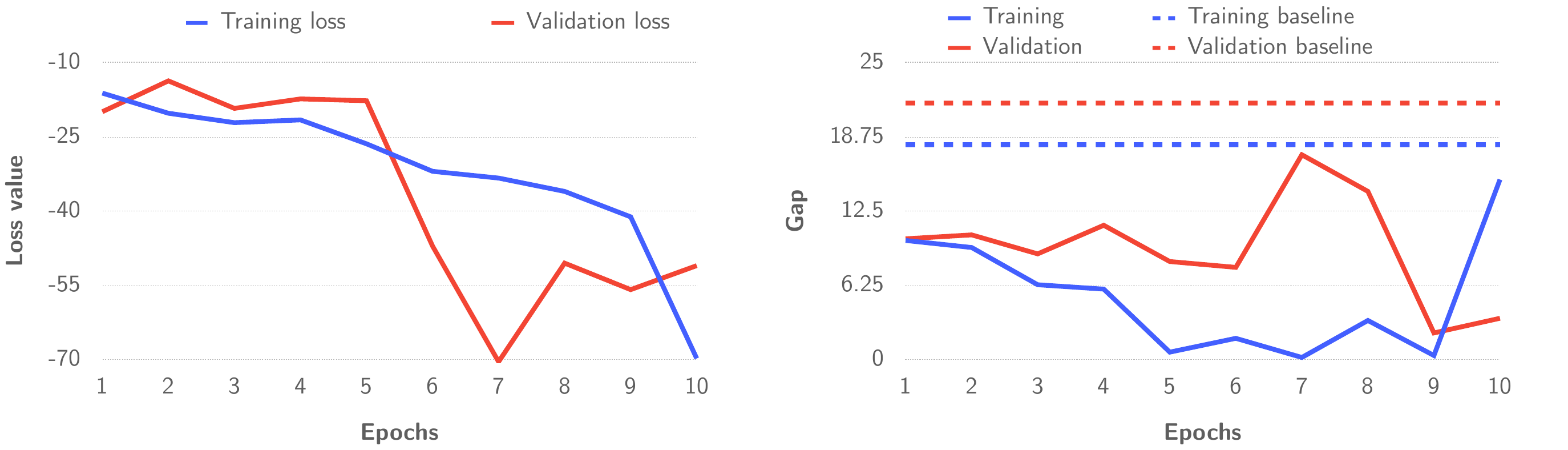}
\caption{Training and validation results for parametric control instances.}
\label{fig:control}
\end{figure}

\paragraph{Results.} The baseline gaps for training, validation, and tests are $18.10\%$, $21.60\%$, and $21.65\%$, respectively. We represent our results in \cref{fig:control,tab:results}. Similarly to 2-matching, our algorithm manages to lower the training gap to a minimum of $0.22\%$ and achieves similar relative improvements on the validations set. We, therefore, stop the training at the $9$th epoch. We observe a training, validation, and test gap of $0.38\%$, $2.28\%$, and $5.47\%$, a remarkable improvement over the baseline. Even in this case, we report an insightful metric on the cut violation: while the algorithm starts with a normalized cut violation of almost $0.45$, it quickly brings it down to an average of $0.25$ after the $5$th epoch. 

\section{Concluding remarks}
\label{sec:theend}
We presented a differentiable cutting-plane algorithm designed to exploit the repetitive structure of parametric mixed-integer optimization. 
Although our computational tests are still preliminary and performed on relatively small sets of instances, we are optimistic about this new research direction.
We believe our work has three areas of potential impact. First, our work can lay the foundations for embedding mixed-integer solvers \emph{as a layer} in differentiable architectures. Second, our algorithm can inform researchers in our community about the structural properties of cutting planes on parametric instances and, perhaps, propel methodological advances in reoptimization. Finally, our approach can potentially enable the solution of small-to-medium sized \glspl{MIP} in embedded devices, such as microcontrollers.

\newpage
\subsubsection*{Acknowledgements} Gabriele Dragotto was supported by the \emph{Schmidt DataX} Fund at Princeton University, made possible through a major gift from the Schmidt Futures Foundation.

\bibliographystyle{splncs04.bst}
\bibliography{bibliography}  %

\begin{thebibliography}{10}
\providecommand{\url}[1]{\texttt{#1}}
\providecommand{\urlprefix}{URL }
\providecommand{\doi}[1]{https://doi.org/#1}

\bibitem{cvxpylayers2019}
Agrawal, A., Amos, B., Barratt, S., Boyd, S., Diamond, S., Kolter, Z.: \href{https://arxiv.org/pdf/1910.12430.pdf}{Differentiable convex optimization layers}. In: Advances in Neural Information Processing Systems (2019)

\bibitem{diffcp2019}
Agrawal, A., Barratt, S., Boyd, S., Busseti, E., Moursi, W.: \href{https://arxiv.org/abs/1904.09043}{Differentiating through a cone program}. Journal of Applied and Numerical Optimization  \textbf{1}(2),  107--115 (2019)

\bibitem{cocp}
Agrawal, A., Barratt, S., Boyd, S., Stellato, B.: \href{https://proceedings.mlr.press/v120/agrawal20a.html}{Learning convex optimization control policies}. In: Proceedings of the 2nd Conference on Learning for Dynamics and Control. Proceedings of Machine Learning Research, vol.~120, pp. 361--373. PMLR (2020)

\bibitem{alvarez2017}
Alvarez, A.M., Louveaux, Q., Wehenkel, L.: \href{https://pubsonline.informs.org/doi/10.1287/ijoc.2016.0723}{A machine learning-based approximation of strong branching}. INFORMS Journal on Computing  \textbf{29}(1),  185--195 (2017)

\bibitem{amaldi2014coordinated}
Amaldi, E., Coniglio, S., Gualandi, S.: \href{https://link.springer.com/article/10.1007/s10107-012-0596-x}{Coordinated cutting plane generation via multi-objective separation}. Mathematical Programming  \textbf{143},  87--110 (2014)

\bibitem{amos_amortized_2023}
Amos, B.: \href{http://dx.doi.org/10.1561/2200000102}{Tutorial on amortized optimization}. Foundations and Trends in Machine Learning  \textbf{16}(5),  592--732 (2023)

\bibitem{balas_disjunctive_2018}
Balas, E.: {Disjunctive Programming}. Springer International Publishing, 1st ed. edn. (2018)

\bibitem{balas_1993_lift}
Balas, E., Ceria, S., Cornu{\'e}jols, G.: \href{https://link.springer.com/content/pdf/10.1007/BF01581273.pdf}{A lift-and-project cutting plane algorithm for mixed 0--1 programs}. Mathematical programming  \textbf{58}(1-3),  295--324 (1993)

\bibitem{balas1996mixed}
Balas, E., Ceria, S., Cornu{\'e}jols, G.: \href{https://pubsonline.informs.org/doi/abs/10.1287/mnsc.42.9.1229}{Mixed 0-1 programming by lift-and-project in a branch-and-cut framework}. Management Science  \textbf{42}(9),  1229--1246 (1996)

\bibitem{balas1980strengthening}
Balas, E., Jeroslow, R.G.: \href{https://www.sciencedirect.com/science/article/abs/pii/037722178090106X}{Strengthening cuts for mixed integer programs}. European Journal of Operational Research  \textbf{4}(4),  224--234 (1980)

\bibitem{balas2008optimizing}
Balas, E., Saxena, A.: \href{https://link.springer.com/article/10.1007/s10107-006-0049-5}{Optimizing over the split closure}. Mathematical Programming  \textbf{113}(2),  219--240 (2008)

\bibitem{MachineLearninBengio2021}
Bengio, Y., Lodi, A., Prouvost, A.: \href{https://www.sciencedirect.com/science/article/abs/pii/S0377221720306895}{Machine learning for combinatorial optimization: A methodological tour d'horizon}. European Journal of Operational Research  \textbf{290},  405--421 (2021)

\bibitem{berthold2022learning}
Berthold, T., Francobaldi, M., Hendel, G.: \href{https://arxiv.org/abs/2206.11618}{Learning to use local cuts}. arXiv  \textbf{2206.11618} (2022)

\bibitem{bertsimas2021}
{Bertsimas}, D., {Stellato}, B.: \href{https://doi.org/10.1007/s10994-020-05893-5}{The voice of optimization}. Machine Learning  \textbf{110},  249--277 (2021)

\bibitem{stellato2022}
{Bertsimas}, D., {Stellato}, B.: \href{https://doi.org/10.1287/ijoc.2022.1181}{Online mixed-integer optimization in milliseconds}. INFORMS Journal on Computing  \textbf{34}(4),  2229--2248 (2022)

\bibitem{scip}
Bestuzheva, K., Besan{\c{c}}on, M., Chen, W.C., Chmiela, A., Donkiewicz, T., van Doornmalen, J., Eifler, L., Gaul, O., Gamrath, G., Gleixner, A., Gottwald, L., Graczyk, C., Halbig, K., Hoen, A., Hojny, C., van~der Hulst, R., Koch, T., L{\"u}bbecke, M., Mahemr, S.J., Matter, F., M{ \"u}hmer, E., M{\"u}ller, B., Pfetsch, M.E., Rehfeldt, D., Schlein, S., Schl{\"o}sser, F., Serrano, F., Shinano, Y., Sofranac, B., Turner, M., Vigerske, S., Wegscheider, F., Wellner, F., Weninger, D., Witzig, J.: \href{http://www.optimization-online.org/DB\_HTML/2021/12/8728.html}{The {SCIP} optimization suite 8.0}. Technical report, Optimization Online (2021)

\bibitem{bixby2010}
Bixby, R.E.: \href{https://www.math.uni-bielefeld.de/documenta/vol-ismp/25_bixby-robert.pdf}{A brief history of linear and mixed-integer programming computation}. Documenta Mathematica pp. 107--121 (2010)

\bibitem{bixby2002solving}
Bixby, R.E.: \href{https://pubsonline.informs.org/doi/abs/10.1287/opre.50.1.3.17780}{Solving real-world linear programs: A decade and more of progress}. Operations research  \textbf{50}(1),  3--15 (2002)

\bibitem{bonami2012optimizing}
Bonami, P.: \href{https://link.springer.com/article/10.1007/s12532-012-0037-0}{On optimizing over lift-and-project closures}. Mathematical Programming Computation  \textbf{4},  151--179 (2012)

\bibitem{optmethodsml}
Bottou, L., Curtis, F.E., Nocedal, J.: \href{https://epubs.siam.org/doi/10.1137/16M1080173}{Optimization methods for large-scale machine learning}. SIAM Review  \textbf{60}(2),  223--311 (Jan 2018)

\bibitem{jax2018github}
Bradbury, J., Frostig, R., Hawkins, P., Johnson, M.J., Leary, C., Maclaurin, D., Necula, G., Paszke, A., Vander{P}las, J., Wanderman-{M}ilne, S., Zhang, Q.: \href{http://github.com/google/jax}{{JAX}: composable transformations of {P}ython+{N}um{P}y programs} (2018)

\bibitem{caprara2003separation}
Caprara, A., Letchford, A.N.: \href{https://link.springer.com/article/10.1007/s10107-002-0320-3}{On the separation of split cuts and related inequalities}. Mathematical Programming  \textbf{94},  279--294 (2003)

\bibitem{cauligi2021coco}
Cauligi, A., Culbertson, P., Schmerling, E., Schwager, M., Stellato, B., Pavone, M.: \href{https://doi.org/10.1109/LRA.2021.3135931}{{CoCo}: Online mixed-integer control via supervised learning}. IEEE Robotics and Automation Letters  \textbf{7}(2),  1447--1454 (2022)

\bibitem{chetelat_continuous_2023}
Ch\'{e}telat, D., Lodi, A.: \href{https://www.sciencedirect.com/science/article/pii/S0167637723001074}{Continuous cutting plane algorithms in integer programming}. Operations Research Letters  \textbf{51}(4),  439--445 (2023)

\bibitem{chvatal1973edmonds}
Chv{\'a}tal, V.: \href{https://www.sciencedirect.com/science/article/pii/0012365X73901672}{Edmonds polytopes and a hierarchy of combinatorial problems}. Discrete mathematics  \textbf{4}(4),  305--337 (1973)

\bibitem{splitcuts}
Cook, W., Kannan, R., Schrijver, A.: \href{https://link.springer.com/article/10.1007/bf01580858}{Chv{\'a}tal closures for mixed integer programming problems}. Mathematical Programming  \textbf{47}(1-3),  155--174 (1990)

\bibitem{cptutorial}
Cornu{\'e}jols, G.: \href{https://link.springer.com/content/pdf/10.1007/s10107-006-0086-0.pdf}{Valid inequalities for mixed integer linear programs}. Mathematical programming  \textbf{112}(1),  3--44 (2008)

\bibitem{cornuejols2002rank}
Cornu{\'e}jols, G., Li, Y.: \href{https://link.springer.com/article/10.1007/s101070100250}{On the rank of mixed 0, 1 polyhedra}. Mathematical Programming  \textbf{91},  391--397 (2002)

\bibitem{dalle2022learning}
Dalle, G., Baty, L., Bouvier, L., Parmentier, A.: \href{https://arxiv.org/pdf/2207.13513.pdf}{Learning with combinatorial optimization layers: a probabilistic approach}. arXiv  \textbf{2207.13513} (2022)

\bibitem{dash2010mir}
Dash, S., G{\"u}nl{\"u}k, O., Lodi, A.: \href{https://link.springer.com/article/10.1007/s10107-008-0225-x}{{MIR} closures of polyhedral sets}. Mathematical Programming  \textbf{121},  33--60 (2010)

\bibitem{dey2018theoretical}
Dey, S.S., Molinaro, M.: \href{https://link.springer.com/article/10.1007/s10107-018-1302-4}{Theoretical challenges towards cutting-plane selection}. Mathematical Programming  \textbf{170},  237--266 (2018)

\bibitem{deza_cutting_2023}
Deza, A., Khalil, E.B.: \href{https://arxiv.org/pdf/2302.09166.pdf}{Machine learning for cutting planes in integer programming: A survey}. arXiv  \textbf{2302.09166} (2023)

\bibitem{fisac2019bridging}
Fern\'{a}ndez~Fisac, J., Lugovoy, N.F., Rubies-Royo, V., Ghosh, S., Tomlin, C.J.: \href{https://ieeexplore.ieee.org/abstract/document/8794107/}{Bridging hamilton-jacobi safety analysis and reinforcement learning}. In: 2019 International Conference on Robotics and Automation (ICRA). pp. 8550--8556. IEEE (2019)

\bibitem{fico}
{FICO, Inc.}: \href{https://www.fico.com/en/products/fico-xpress-optimization}{Optimizer reference manual} (2023)

\bibitem{ontheseparation}
Fischetti, M., Lodi, A., Tramontani, A.: \href{https://link.springer.com/content/pdf/10.1007/s10107-009-0300-y.pdf}{On the separation of disjunctive cuts}. Mathematical Programming  \textbf{128}(1-2),  205--230 (2011)

\bibitem{fischetti2013approximating}
Fischetti, M., Salvagnin, D.: \href{https://pubsonline.informs.org/doi/abs/10.1287/ijoc.1120.0543}{Approximating the split closure}. INFORMS Journal on Computing  \textbf{25}(4),  808--819 (2013)

\bibitem{frick2015embedded}
Frick, D., Domahidi, A., Morari, M.: \href{https://www.sciencedirect.com/science/article/abs/pii/S0098135414001847}{Embedded optimization for mixed logical dynamical systems}. Computers \& Chemical Engineering  \textbf{72},  21--33 (2015)

\bibitem{gamrath2015reoptimization}
Gamrath, G., Hiller, B., Witzig, J.: \href{https://link.springer.com/chapter/10.1007/978-3-319-20086-6_14}{Reoptimization techniques for {MIP} solvers}. In: Experimental Algorithms: 14th International Symposium, SEA 2015, Paris, France, June 29--July 1, 2015, Proceedings 14. pp. 181--192 (2015)

\bibitem{exactcombopt}
Gasse, M., Chetelat, D., Ferroni, N., Charlin, L., Lodi, A.: \href{https://proceedings.neurips.cc/paper/2019/file/d14c2267d848abeb81fd590f371d39bd-Paper.pdf}{Exact combinatorial optimization with graph convolutional neural networks}. In: Advances in Neural Information Processing Systems. vol.~32 (2019)

\bibitem{mipcc23}
Gleixner, A., Berthold, T., Besan\c{c}on, M., Bolusani, S., D'Ambrosio, C., Mu\~{n}oz, G., Paat, J., Thomopulos, D.: \href{https://www.mixedinteger.org/2023}{{MIP} computational competition 2023} (2023)

\bibitem{gomory_outline_1958}
Gomory, R.E.: \href{https://open.spotify.com/playlist/4MmL29PYtR9lm9h9sx7SX8?si=bdfcccc4d05044ca}{Outline of an algorithm for integer solutions to linear programs}. Bull. Amer. Math. Soc.  \textbf{64}(5),  275--279 (1958)

\bibitem{clarabel}
Goulart, P., Chen, Y.: \href{https://oxfordcontrol.github.io/ClarabelDocs}{Clarabel documentation} (2023)

\bibitem{gurobi}
{Gurobi Optimization, LLC}: \href{https://www.gurobi.com}{Gurobi optimizer reference manual} (2023)

\bibitem{guzelsoy2007duality}
Guzelsoy, M., Ralphs, T.K.: \href{https://coral.ise.lehigh.edu/~ted/files/papers/MILPD06.pdf}{Duality for mixed-integer linear programs}. International Journal of Operations Research  \textbf{4}(3),  118--137 (2007)

\bibitem{huang2022learning}
Huang, Z., Wang, K., Liu, F., Zhen, H.L., Zhang, W., Yuan, M., Hao, J., Yu, Y., Wang, J.: \href{https://www.sciencedirect.com/science/article/abs/pii/S0031320321005331}{Learning to select cuts for efficient mixed-integer programming}. Pattern Recognition  \textbf{123},  108353 (2022)

\bibitem{highs}
Huangfu, Q., Hall, J.A.: \href{https://link.springer.com/article/10.1007/s12532-017-0130-5}{Parallelizing the dual revised simplex method}. Mathematical Programming Computation  \textbf{10}(1),  119--142 (2018)

\bibitem{cplex}
{IBM, Inc.}: \href{https://www.ibm.com/products/ilog-cplex-optimization-studio}{Optimizer reference manual} (2023)

\bibitem{kelley1960cutting}
Kelley, Jr, J.E.: \href{https://epubs.siam.org/doi/abs/10.1137/0108053?journalCode=smjmap.1}{The cutting-plane method for solving convex programs}. Journal of the society for Industrial and Applied Mathematics  \textbf{8}(4),  703--712 (1960)

\bibitem{learn2branch}
Khalil, E.B., Bodic, P., Song, L., Nemhauser, G., Dilkina, B.: \href{https://ojs.aaai.org/index.php/AAAI/article/view/10080}{Learning to branch in mixed integer programming}. In: Proceedings of the Thirtieth AAAI Conference on Artificial Intelligence. pp. 724--731. AAAI'16, AAAI Press (2016)

\bibitem{kilincc2016minimal}
K{\i}l{\i}n{\c{c}}-Karzan, F.: \href{https://pubsonline.informs.org/doi/abs/10.1287/moor.2015.0737}{On minimal valid inequalities for mixed integer conic programs}. Mathematics of Operations Research  \textbf{41}(2),  477--510 (2016)

\bibitem{land_automatic_1960}
Land, A.H., Doig, A.G.: \href{https://www.jstor.org/stable/1910129?origin=crossref}{An automatic method of solving discrete programming problems}. Econometrica  \textbf{28}(3), ~497 (1960)

\bibitem{lodi2023disjunctive}
Lodi, A., Tanneau, M., Vielma, J.P.: \href{https://link.springer.com/article/10.1007/s10107-022-01844-1}{Disjunctive cuts in mixed-integer conic optimization}. Mathematical Programming  \textbf{199}(1-2),  671--719 (2023)

\bibitem{Lodi2017}
Lodi, A., Zarpellon, G.: \href{https://link.springer.com/article/10.1007/s11750-017-0451-6}{On learning and branching: a survey}. TOP  \textbf{25}(2),  207--236 (2017)

\bibitem{ralphs_selecting_2010}
Mahajan, A., Ralphs, T.: \href{https://doi.org/10.1137/080737587}{On the complexity of selecting disjunctions in integer programming}. SIAM Journal on Optimization  \textbf{20}(5),  2181--2198 (2010)

\bibitem{warmstart}
Marcucci, T., Tedrake, R.: \href{https://ieeexplore.ieee.org/document/9134792}{Warm start of mixed-integer programs for model predictive control of hybrid systems}. IEEE Transactions on Automatic Control  \textbf{66}(6),  2433--2448 (2021)

\bibitem{moerland2023model}
Moerland, T.M., Broekens, J., Plaat, A., Jonker, C.M.: \href{https://ieeexplore.ieee.org/document/10007800}{Model-based reinforcement learning: A survey}. Foundations and Trends in Machine Learning  \textbf{16}(1),  1--118 (2023)

\bibitem{nemhauser1990recursive}
Nemhauser, G.L., Wolsey, L.A.: \href{https://link.springer.com/article/10.1007/BF01585752}{A recursive procedure to generate all cuts for 0--1 mixed integer programs}. Mathematical Programming  \textbf{46}(1-3),  379--390 (1990)

\bibitem{padberg2013linear}
Padberg, M.: {Linear Optimization and Extensions}, vol.~12. Springer Science \& Business Media (2013)

\bibitem{padberg_branch-and-cut_1991}
Padberg, M., Rinaldi, G.: \href{http://www.jstor.org/stable/2030652}{A branch-and-cut algorithm for the resolution of large-scale symmetric traveling salesman problems}. SIAM Review  \textbf{33}(1),  60--100 (1991)

\bibitem{patel2023progressively}
Patel, K.K.: \href{https://arxiv.org/pdf/2308.08986.pdf}{Progressively strengthening and tuning {MIP} solvers for reoptimization}. arXiv  \textbf{2308.08986} (2023)

\bibitem{paulus2022learning}
Paulus, M.B., Zarpellon, G., Krause, A., Charlin, L., Maddison, C.: \href{https://proceedings.mlr.press/v162/paulus22a/paulus22a.pdf}{Learning to cut by looking ahead: Cutting plane selection via imitation learning}. In: International conference on machine learning. pp. 17584--17600 (2022)

\bibitem{sutton2018reinforcement}
Sutton, R.S., Barto, A.G.: Reinforcement learning: An introduction. MIT press (2018)

\bibitem{faenza_cutting_2020}
Tang, Y., Agrawal, S., Faenza, Y.: \href{https://proceedings.mlr.press/v119/tang20a.html}{Reinforcement learning for integer programming: Learning to cut}. In: III, H.D., Singh, A. (eds.) Proceedings of the 37th International Conference on Machine Learning. Proceedings of Machine Learning Research, vol.~119, pp. 9367--9376. PMLR (2020)

\bibitem{turner2023cutting}
Turner, M., Berthold, T., Besan\c{c}on, M., Koch, T.: \href{https://link.springer.com/chapter/10.1007/978-3-031-33271-5_4}{Cutting plane selection with analytic centers and multiregression}. In: International Conference on Integration of Constraint Programming, Artificial Intelligence, and Operations Research. pp. 52--68 (2023)

\bibitem{turner2022adaptive}
Turner, M., Koch, T., Serrano, F., Winkler, M.: \href{https://arxiv.org/pdf/2202.10962.pdf}{Adaptive cut selection in mixed-integer linear programming}. arXiv  \textbf{2202.10962} (2022)

\bibitem{walter2013sparsity}
Walter, M.: \href{https://link.springer.com/chapter/10.1007/978-3-319-00795-3_2}{Sparsity of lift-and-project cutting planes}. In: Operations Research Proceedings 2012: Selected Papers of the International Annual Conference of the German Operations Research Society (GOR), Leibniz University of Hannover, Germany, September 5-7, 2012. pp. 9--14. Springer (2013)

\bibitem{werbos1990backpropagation}
Werbos, P.: \href{https://ieeexplore.ieee.org/document/58337}{Backpropagation through time: {W}hat it does and how to do it}. Proc. IEEE  \textbf{78}(10),  1550--1560 (1990)

\bibitem{reinforce}
Williams, R.J.: \href{https://link.springer.com/content/pdf/10.1007/BF00992696.pdf}{Simple statistical gradient-following algorithms for connectionist reinforcement learning}. Machine learning  \textbf{8},  229--256 (1992)

\end{thebibliography}

\newpage
\appendix

\section{Appendix}

\subsection{Relationship to subadditive cutting planes}

\begin{proof}[\cref{thm:continuouscp}]
\noindent \emph{Proof of \ref{thm:continuouscp:1}. } This statement follows from the definition of \cref{eq:CGP} and the theorem's assumptions. \\
\noindent \emph{Proof of \ref{thm:continuouscp:2}. }
Let $\gomorydisjpart_j = -\doubleu^\transpose \probA_j$. By substituting the appropriate terms from the definition, we have that
\begin{align*}
    \cutg_j &= 
        \min \Bigg(\doubleuplus^\transpose \begin{bmatrix}
        \probA_j \\ \disj_j
    \end{bmatrix}, \doubleuminus^\transpose \begin{bmatrix}
        \probA_j \\ -\disj_j
    \end{bmatrix}\Bigg) \\
    &= \max \Big( \doubleuplus^\transpose \probA_j + \uzero \lfloor \gomorydisjpart_j \rfloor, \doubleuminus^\transpose \probA_j - \vzero \lceil \gomorydisjpart_j \rceil \Big).
\end{align*}
By combining the fact that $\doubleuplus^\transpose \probA_j = \uzero (\doubleuplus - \doubleuminus)^\transpose \probA_j + (\vzero \doubleuplus + \uzero \doubleuminus)^\transpose \probA_j$, and the fact that $\doubleuminus^\transpose \probA_j = - \vzero (\doubleuplus - \doubleuminus)^\transpose \probA_j + (\vzero \doubleuplus + \uzero \doubleuminus)^\transpose \probA_j$ with the definition of $\gomorydisjpart$, we see that,
\begin{align*}
    \cutg_j &= \max \Big(\uzero (-\gomorydisjpart_j + \lfloor \gomorydisjpart_j \rfloor) + (\vzero \doubleuplus + \uzero \doubleuminus)^\transpose\probA_j, \\
    & \qquad \qquad 
    \vzero (- \lceil \gomorydisjpart_j \rceil + \gomorydisjpart_j) + (\vzero \doubleuplus + \uzero \doubleuminus)^\transpose\probA_j \Big) \\
    &= \max \Big(\uzero (-\gomorydisjpart_j + \lfloor \gomorydisjpart_j \rfloor),
    \vzero (- \lceil \gomorydisjpart_j \rceil + \gomorydisjpart_j) \Big)
     + (\vzero \doubleuplus + \uzero \doubleuminus)^\transpose\probA_j.
\end{align*}
Since $\doubleuplus$ is complementary to $\doubleuminus$, then 
$$\vzero \doubleuplus^\transpose + \uzero \doubleuminus^\transpose = \max(\vzero \doubleuplus, \uzero \doubleuminus) = \max((1 - \ve) \doubleu, -\ve \doubleu)^\transpose.$$
Using this latter fact and the definition of $\gomorydisjpart$, then
\begin{align*}
     \cutg_j
    &= \max \Big((1 - \ve) (\doubleu \probA_j + \lfloor -\doubleu \probA_j \rfloor), 
    \ve (-\lceil -\doubleu \probA_j \rceil - \doubleu \probA_j)\Big)
      \\
      &\qquad \qquad  + \max\Big(\ve \doubleu, -(1 - \ve)\doubleu\Big)^\transpose \probA_j \\
    &= \max \Big((1 - \ve) \big(-\fractionalpart{- \doubleu^\transpose \probA_j} \big), \\
    &\qquad \qquad \ve \big(\fractionalpart{ - \doubleu^\transpose \probA_j} - 1\big)\Big)
       \max\Big(\ve \doubleu, -(1 - \ve) \doubleu\Big)^\transpose \probA_j \\
    &= - \min \Big(\ve \fractionalpart{- \doubleu^\transpose \probA_j},
    (1 - \ve) \big(1 - \fractionalpart{ - \doubleu^\transpose \probA_j}\big)\Big) +\\
    & \qquad \qquad
      - \max\Big(\ve \doubleu, -(1 - \ve) \doubleu\Big)^\transpose (-\probA_j).
\end{align*}
If we substitute $- \probA_j$ into the definition of $\subaddfn_{\doubleu, \ve}$, we obtain
\begin{align*}
    \cutg_j &= - \subaddfn_{\doubleu, \ve}(-\probA_j)(1 - \ve).
\end{align*}
We now turn our attention to the right-hand side. Using the definition of $\cuth, \vzero$ and $\disjvalue$, we see that $\probb^\transpose \doubleu + \vzero = \disjvalue$, or equivalently, $\probb^\transpose \doubleuplus + \uzero \disjvalue = \probb^\transpose \doubleuminus - \vzero(\disjvalue + 1)$. Using this equality and the definition of $\cuth$,
\begin{align*}
    \cuth &=  \begin{bmatrix}\probb^\transpose  & -(\disjvalue+1) \end{bmatrix}\cutv \\
    &= \uzero \begin{bmatrix}\probb^\transpose  & -(\disjvalue+1) \end{bmatrix}\cutv + \vzero \begin{bmatrix}\probb^\transpose  & -(\disjvalue+1) \end{bmatrix}\cutv\\
    &= \uzero \begin{bmatrix}\probb^\transpose  & -(\disjvalue+1) \end{bmatrix}\cutv + \vzero \begin{bmatrix} \probb^\transpose & \disjvalue  \end{bmatrix}\cutu.
\end{align*}
We now employ the definition of $\cutu$ and $\cutv$, and obtain
\begin{align*}
    \cuth &= \uzero (\probb^\transpose \doubleuminus - (\disjvalue + 1)\vzero ) +     \vzero(\probb^\transpose \doubleuplus + \disjvalue \uzero ) \\
    &= \uzero \probb^\transpose \doubleuminus + \vzero \probb^\transpose \doubleuplus + \uzero \vzero.
\end{align*}
Similarly to the left-hand side, we note that $\uzero \doubleuminus + \vzero \doubleuplus = \max(-\ve \doubleu, (1-\ve)\doubleu)$. We substitute  $\uzero$ and $\vzero$, and use the definition of $\subaddfn_{\doubleu, \ve}$ to get
\begin{align*}
    \cuth &= \max(-\ve \doubleu, (1-\ve)\doubleu)^\transpose \probb - \ve(1-\ve) \\
    &= - \subaddfn_{\doubleu, \ve}(-\probb)(1 - \ve). 
\end{align*}
Therefore, the inequality $\cutg^\transpose \variables \le \cuth$ is equivalent to $- (1 - \ve) \subaddfn_{\doubleu, \ve}(-\probA) \variables \le - (1 - \ve)\subaddfn_{\doubleu, \ve}(-\probb)$. 
\end{proof}

\subsection{Norm-2 representability}
\begin{proof}[\cref{coro:norm2}]
    Because of \cref{thm:continuouscp}, every subadditive cut is in the feasible region of \cref{eq:CGP}. We need to show that with $\normtype=2$ and some choice of $(\normalization, \disj, \disjvalue)$, any feasible point of \cref{eq:CGP} with a positive objective value can become optimal. Let $\doubleu$ and $\ve$ be such that the subadditive cut $-\subaddfn_{\doubleu, \ve}(-\probA) \variables \le -\subaddfn_{\doubleu, \ve}(-\probb)$ cuts off $\toseparate$, \ie, such that $-\subaddfn_{\doubleu, \ve}(-\probA) \toseparate > -\subaddfn_{\doubleu, \ve}(-\probb)$. Fix $\disj, \disjvalue$ as defined in \cref{thm:continuouscp} and let $\cutu^\star = (\doubleuplus,  \uzero), \cutv^\star=(\doubleuplus,  \vzero)$, where, as before, $\vzero = \fractionalpart{\doubleu^\transpose \probb}$ and $\uzero = 1 - \vzero$. Let
    \begin{equation*}
        \spareMone =\begin{bmatrix}
            \probA \toseparate \\ \disj^\transpose \toseparate
        \end{bmatrix},
        \quad \spareMtwo =\begin{bmatrix}
            \probA \toseparate \\ -\disj^\transpose \toseparate
        \end{bmatrix}, \quad
        \spareMthree = - \begin{bmatrix}
            \probb \\ \disjvalue \end{bmatrix}, 
         \quad \spareMfour = - \begin{bmatrix}
             \probb \\ (\disjvalue + 1) \end{bmatrix}.
    \end{equation*}
   We can rewrite the \gls{CGP} \cref{eq:CGP} as
    \begin{equation} \label{eq:CGP2}
	\begin{array}{ll}
        \underset{{\cutu,\cutv}}{\maxi} \quad & \min \big( \spareMone^\transpose \cutu + \spareMthree^\transpose \cutu, \spareMtwo^\transpose \cutv + \spareMfour^\transpose \cutv, \spareMone^\transpose \cutu + \spareMfour^\transpose \cutv, \spareMtwo^\transpose \cutu + \spareMthree^\transpose \cutv \big) \\[.5em]
		\st \quad & \quad \left\lVert \normalization\begin{bmatrix} \cutu \\ \cutv \end{bmatrix}\right\lVert_2\le1, \\
        & \cutu \ge 0, \cutv \ge 0.
	\end{array}
    \end{equation}
    Without loss of generality, we assume $\spareMone^\transpose \cutu^\star + \spareMthree^\transpose \cutu^\star$ attains the minimum in \cref{eq:CGP2}, and, therefore, that $\spareMone^\transpose \cutu^\star + \spareMthree^\transpose \cutu^\star > 0$.
    We need to choose $\normalization$ such that $\cutu^\star, \cutv^\star$ are optimal in \cref{eq:CGP2}, or equivalently such that $(\cutu^\star, \cutv^\star)$ is optimal in,
    \begin{equation*}
    \label{eq:CGP3}
	\begin{array}{ll}
        \underset{{\cutu,\cutv}}{\maxi} & \spareMone^\transpose \cutu + \spareMthree^\transpose \cutu \\
		\st  &  \left\lVert \normalization\begin{bmatrix} \cutu \\ \cutv \end{bmatrix}\right\lVert_2^2\le1\\
        & \cutu \ge 0, \cutv \ge 0.
	\end{array}
    \end{equation*}
    Let $\utilde = (\cutu^\star, \cutv^\star)$ and $\spareMfive = (\spareMone + \spareMthree, \zerovec_\numconstraints)$, where $\zerovec_\numconstraints$ is a vector of $\numconstraints$ zeroes. If we can choose $\normalization$ such that, for some $\lambda > 0$,
    \begin{align*}
        2 \normalization^\transpose \normalization \utilde = \lambda \spareMfive, \quad \left\lVert \normalization \utilde \right\lVert_2^2 = 1,
    \end{align*}
    holds, we completed our proof. Let 
    $$\tilde{\normalization} = \frac{1}{2\utilde^\transpose \spareMfive} \spareMfive \spareMfive^\transpose.$$ 
    Then, $\lambda\tilde{\normalization}$ is positive semi-definite and $2 \lambda \tilde{\normalization} \utilde = \lambda \spareMfive$ for any $\lambda > 0$. Therefore, we can take $\normalization$ to be the matrix square-root of $\lambda \tilde{\normalization}$ and choose $\lambda>0$ such that $\left\lVert \normalization \utilde \right\lVert_2^2 = 1$.
    By the Karush–Kuhn–Tucker optimality conditions, this completes the proof.
\end{proof}

\subsection{Optimal control problem formulation}
\label{app:control}
We aim to control the energy balance between a supercapacitor and a fuel cell to match a given power demand. The goal is to minimize energy losses while maintaining the device within
acceptable operating limits to prevent lifespan degradation. For a detailed description of the problem, we refer the reader to~\cite{stellato2022,frick2015embedded}. 
The optimization problem reads as
\begin{equation*}
\label{eq:controlproblem}
	\begin{array}{ll}
        \underset{\power_\controltime, \fuelstate_\controltime, \dcontrol_\controltime, \wcontrol_\controltime}{\mini\;} \quad &\sum_{\controltime=0}^{\controlmaxtime-1} \controlf(\power_\controltime, \fuelstate_\controltime)\\
		\st \quad & \energy_{\controltime + 1} = \energy_\controltime + \samplingtime (\power_\controltime - \powerload_{\controltime}) \\
    & \minenergy \leq \energy_{\controltime} \leq \maxenergy \\
    & 0 \leq \power_\controltime \leq \fuelstate_\controltime \maxpower \\
    & \fuelstate_{\controltime+1} = \fuelstate_\controltime + \wcontrol_\controltime \\[-3pt]
    & \scontrol_{\controltime+1} = \scontrol_\controltime + \dcontrol_\controltime - \dcontrol_{\controltime-\controlmaxtime} \\
    & \scontrol_\controltime \leq \nsw \\
    & \Gcontrol(\wcontrol_\controltime, \fuelstate_\controltime, \dcontrol_\controltime) \leq \hcontrol \\
    & \energy_0 = \initenergy, \quad \fuelstate_0 = \fuelstateinit, \quad \scontrol_0 = \scontrolinit, \\
    & \fuelstate \in \{0,1\}^\controlmaxtime, \quad \dcontrol \in \{0, 1\}^\controlmaxtime, \quad \wcontrol \in [-1, 1]^\controlmaxtime.
	\end{array}
\end{equation*}
Here, $\controltime$ > 0 is the sampling time, $\controlmaxtime$ is the horizon length, $\energy_\controltime \in [\minenergy, \maxenergy]$ is the energy stored, $\power_\controltime \in [0, \maxpower]$ is the power provided by the fuel cell, $\powerload_\controltime$ is the desired power load. 
The variable $\fuelstate_\controltime \in \{0,1\}$ is the on-off state of the fuel cell, whereas 
$\dcontrol_\controltime \in \{0, 1\}$ is a binary variable which is $1$ iff the control $\fuelstate_\controltime$ changes at time $\controltime$ 
and $0$ otherwise.
The variable $\wcontrol_\controltime$ represents the amount of change brought by 
$\dcontrol_\controltime$ in the battery state 
$\fuelstate_\controltime$, while $\scontrol_\controltime$ is the number of times the control switched between $\controltime - \controlmaxtime$ and $\controltime$. We set a limit $\nsw$ to the maximum allowed number of switches between times $\controltime - \controlmaxtime$ and $\controltime$. The problem parameters are $\parameters = (\initenergy,
\fuelstateinit, \scontrolinit, \dpast, 
\powerload)$ where $\dpast = (\dcontrol_{-\controlmaxtime}, \dots, \dcontrol_{-1})$ and $\powerload = (\powerload_0, \dots \powerload_{\controlmaxtime - 1})$. The objective term
$\controlf$ is a function representing the affine stage power cost, \ie, $\controlf(\power, \fuelstate) = \betacontrol \power + \gammacontrol \fuelstate$ for fixed parameters $\betacontrol, \gammacontrol \in \R$. 
\end{document}